\documentclass[a4paper,12pt]{article}
\usepackage{amsthm,amssymb,amsmath,fullpage,url}

\newtheorem{theorem}{Theorem}
\newtheorem{lemma}{Lemma}
\newtheorem{proposition}{Proposition}
\newtheorem{corollary}{Corollary}
\newtheorem{conjecture}{Conjecture}
\theoremstyle{remark}
\newtheorem{remark}{Remark}
\theoremstyle{definition}
\newtheorem{definition}{Definition}

\title{Pattern Occurrence in the Dyadic Expansion of Square Root of Two and an Analysis of Pseudorandom Number Generators}
\author{Koji Nuida}
\date{Research Center for Information Security (RCIS), National Institute of Advanced Industrial Science and Technology (AIST), Japan\\\url{k.nuida@aist.go.jp}}
\begin{document}
\maketitle
\begin{abstract}
Recently, designs of pseudorandom number generators (PRNGs) using integer-valued variants of logistic maps and their applications to some cryptographic schemes have been studied, due mostly to their ease of implementation and performance.
However, it has been noted that this ease is reduced for some choices of the PRNGs accuracy parameters.
In this article, we show that the distribution of such undesirable accuracy parameters is closely related to the occurrence of some patterns in the dyadic expansion of the square root of 2.
We prove that for an arbitrary infinite binary word, the asymptotic occurrence rate of these patterns is bounded in terms of the asymptotic occurrence rate of zeroes.
We also present examples of infinite binary words that tightly achieve the bounds.
As a consequence, a classical conjecture on asymptotic evenness of occurrence of zeroes and ones in the dyadic expansion of the square root of 2 implies that the asymptotic rate of the undesirable accuracy parameters for the PRNGs is at least 1/6.
\end{abstract}

\section{Introduction}
\label{sec:intro}

Randomness is a ubiquitous element in our present life and can be found from a simple coin toss at the beginning of a football game, to more complex settings such as encrypted communication of governmental secrets.
The provision, application and evaluation of randomness has occupied a major and attractive branch of mathematics.
In particular, there exist several methods and techniques that generate a seemingly random-looking sequence by using shorter random sequences (often called a \emph{seed}) and deterministic algorithm, better known as a \emph{pseudorandom number generator} (\emph{PRNG}).
See, for instance, \cite{Knu_book} for references.

In this article, we reveal a nontrivial relation between analysis of some PRNGs and properties of the dyadic expansion of $\sqrt{2} = (1.01101 \cdots)_2$.
Note, however, that the dyadic expansion of $\sqrt{2}$ does not appear in the construction of the PRNGs itself.
Also, it has been shown that the logistic map
\begin{displaymath}
L(x) = \mu x (1-x) \enspace,\enspace 0 < x < 1 \enspace,
\end{displaymath}
for some parameter $\mu$ can be effectively used for constructing good PRNGs (\cite{PR95,Wag93}).
In particular, when $\mu = 4$ is adopted, the logistic map shows chaotic behavior.
However, those PRNGs deal with real number as outputs and can therefore not be implemented in computers due to their finite accuracy.
As a result, a modified integer-valued logistic map of the form:
\begin{displaymath}
L_n(x) = \left\lfloor \frac{ 4 x (2^n - x) }{ 2^n } \right\rfloor = \left\lfloor \frac{ x (2^n - x) }{ 2^{n-2} } \right\rfloor \enspace,\enspace x \in X_n = \{1,2,\dots,2^n - 1\}
\end{displaymath}
where $\lfloor z \rfloor$ denotes the largest integer $N$ such that $N \leq z$ and $2 \leq n \in \mathbb{Z}$ is an accuracy parameter, has been proposed and studied in \cite{AMU06,MAU08}.
The definition of $L_n(x)$ is derived from $L(x)$ by expanding the bounds of the original seed $x \in (0,1)$ to the larger interval $(0,2^n)$ and then truncating the final value to obtain an integer.
The corresponding PRNG first chooses an internal state $s_0 = s$ from the set $X_n$ and then for each step $i \geq 1$, updates the internal state by $s_i = L_n(s_{i-1})$ and outputs some bits in the dyadic expansion of $s_i$.

For the above PRNG, it has been mentioned in \cite{MAU07} that when $s_i = 2^{n-1}$ for some $i$, the subsequent internal states become eventually stable, namely we have $s_{i+1} = 2^n$ and $s_j = 0$ for every $j \geq i + 2$.
Since internal states of PRNGs being stable are fatal for the purpose of providing good randomness, the value $2^{n-1}$ should not be used as an internal state.
To concern the problem, in general it is not sufficient to exclude the value $2^{n-1}$ itself from the candidates of the initial internal state $s_0$.
Namely, if there exists an $x \in X_n$ such that $L_n(x) = 2^{n-1}$, then the choice of internal state $s_0 = x$ for such an $x$ also makes the internal states eventually stable.
We call the accuracy parameter $n$ \emph{undesirable} if such an $x$ exists, since in such a case an extra check is required for choosing an appropriate initial internal state.
The motivation of this work is to estimate how many undesirable parameters exist among the integers $n \geq 2$.

We explain an aforementioned relation of the above PRNGs with combinatorial properties of $\sqrt{2}$.
Let $b_i \in \{0,1\}$ denote the $i$-th bit of the fractional part of the dyadic expansion of $\sqrt{2}$, namely
\begin{displaymath}
\sqrt{2} = (1.b_1 b_2 b_3 \cdots)_2 \enspace.
\end{displaymath}
We show that a parameter $n \geq 2$ is undesirable if the ($n-1$)-th tail $b_{n-1} b_n b_{n+1} \cdots$ of the dyadic expansion of $\sqrt{2}$ begins with one of the three patterns $00$, $0100$, and $01010$.
For instance, since
\begin{equation}
\label{eq:dyadic_expansion_sqrt2}
\sqrt{2} = (1.011010100000100 \cdots)_2 \enspace,
\end{equation}
we have $b_{12} = b_{14} = b_{15} = 0$ and $b_{13} = 1$, therefore the above fact implies that $n = 13$ is undesirable.
As a result, the occurrence rate of these three patterns in the dyadic expansion of $\sqrt{2}$ gives a lower bound of the occurrence rate of undesirable parameters.
Motivated by the observation, we study the distributions of the three patterns in \emph{arbitrary} infinite binary words $w = w_1 w_2 w_3 \cdots$, and prove that the asymptotic occurrence rate of the three patterns in $w$ is bounded by a function of the asymptotic occurrence rate of zeroes in $w$.
We also present construction of infinite binary words that achieve the bounds tightly.
(See Theorem \ref{thm:main} for the precise statement.)
This result connects the asymptotic occurrence rate of undesirable parameters to the distribution of zeroes in the dyadic expansion of $\sqrt{2}$.
For the latter, it has been conjectured that the asymptotic occurrence rate of zeroes in the dyadic expansion of $\sqrt{2}$ is $1/2$ (in other words, $\sqrt{2}$ is simply normal to the base $2$).
If the conjecture is true, it follows that the asymptotic occurrence rate of undesirable parameters is lower bounded by $1/6$, which shows a disadvantage of the above PRNGs.

This article is organized as follows.
In Section \ref{sec:logistic}, we prove the aforementioned sufficient condition of an accuracy parameter $n$ being undesirable, in terms of the occurrence rate of the three patterns $00$, $0100$, and $01010$ in the dyadic expansion of $\sqrt{2}$.
In Section \ref{sec:occurrence}, we state the main theorem (Theorem \ref{thm:main}) on a relation between the asymptotic occurrence rates of the three patterns and of zeroes in arbitrary infinite binary words.
As a result, we also estimate the asymptotic occurrence rate of undesirable parameters.
Finally, Section \ref{sec:proof} gives the proof of the main theorem.
\paragraph*{Acknowledgments.}
The author would like to thank Dr.\ Yoshio Okamoto, Dr.\ Kenji Kashiwabara, and Professor Masahiro Hachimori, for their significant comments.
An extended abstract of this work has been accepted by 21st International Conference on Formal Power Series and Algebraic Combinatorics (FPSAC 2009) \cite{Nui_FPSAC}.

\section{Integer-Valued Logistic Maps}
\label{sec:logistic}

As we have mentioned in Section \ref{sec:intro}, a main object of this article is integer-valued logistic maps $L_n(x)$ with domain $X_n = \{1,2,\dots,2^n-1\}$, parameterized by an integer $n \geq 2$, defined by
\begin{equation}
\label{eq:logistic_map}
L_n(x) = \left\lfloor \frac{ 4 x (2^n - x) }{ 2^n } \right\rfloor = \left\lfloor \frac{ x (2^n - x) }{ 2^{n-2} } \right\rfloor \enspace,\enspace x \in X_n = \{1,2,\dots,2^n - 1\} \enspace.
\end{equation}
Note that $L_n(x) \in X_n$ for any $x \in X_n \setminus \{2^{n-1}\}$, while $L_n(2^{n-1}) = 2^n$.
We would like to estimate the asymptotic occurrence rate of accuracy parameters $n$, among all integers $n \geq 2$, that satisfy the following condition:
\begin{definition}
\label{defn:undesirable}
We say that a parameter $2 \leq n \in \mathbb{Z}$ is \emph{undesirable} if there exists an $x \in X_n$ such that $L_n(x) = 2^{n-1}$.
\end{definition}
This definition is motivated by an analysis of some pseudorandom number generators (PRNGs) using $L_n(x)$; see Section \ref{sec:intro} for details.
In the rest of this section, we show that the occurrence rate of undesirable parameters is related to the occurrence of the patterns $00$, $0100$, and $01010$ in the dyadic expansion of $\sqrt{2}$.
For the purpose, first note that by the definition, a parameter $n$ is undesirable if and only if there exists an $x \in X_n$ such that $2^{n-1} \leq x (2^n - x) / 2^{n-2} < 2^{n-1} + 1$.
By solving the inequality, it follows that this condition for $x$ is equivalent to
\begin{equation}
\label{eq:dangerous_condition_1}
\sqrt{ 2^{2n-3} - 2^{n-2} } < |2^{n-1} - x| \leq \sqrt{ 2^{2n-3} } \enspace.
\end{equation}
Moreover, since
\begin{displaymath}
\sqrt{ 2^{2n-3} } - \sqrt{ 2^{2n-3} - 2^{n-2} }
= \frac{ 2^{n-2} }{ \sqrt{ 2^{2n-3} } + \sqrt{ 2^{2n-3} - 2^{n-2} } }
> \frac{ 2^{n-2} }{ 2 \sqrt{ 2^{2n-3} } } = \frac{ \sqrt{2} }{ 4 } \enspace,
\end{displaymath}
the condition (\ref{eq:dangerous_condition_1}) is satisfied if $2^{n-2}\sqrt{2} - \sqrt{2}/4 \leq |2^{n-1} - x| \leq 2^{n-2}\sqrt{2}$.
Summarizing, we have the following lemma:
\begin{lemma}
\label{lem:dangerous_sufficient}
A parameter $n \geq 2$ is undesirable if $2^{n-2}\sqrt{2} - \sqrt{2} / 4 \leq m \leq 2^{n-2}\sqrt{2}$ for some integer $m$.
\end{lemma}
This lemma can be rephrased in terms of the dyadic expansion of $\sqrt{2}$ as follows.
Let $\sqrt{2} = (1.b_1 b_2 b_3 \cdots)_2$ be the dyadic expansion of $\sqrt{2}$.
For instance, we have $b_1 = 0$, $b_2 = 1$ and $b_3 = 1$ (see (\ref{eq:dyadic_expansion_sqrt2})).
Then the fractional part of the dyadic expansion of $2^{n-2} \sqrt{2}$ is $(0.b_{n-1}b_n b_{n+1} \cdots)_2$, while the dyadic expansion of $\sqrt{2} / 4$ is $(0.01b_1 b_2 b_3\cdots)_2$.
By using these expressions, Lemma \ref{lem:dangerous_sufficient} implies the following:
\begin{lemma}
\label{lem:dangerous_sufficient_dyadic}
In the above setting, a parameter $n \geq 2$ is undesirable if
\begin{equation}
\label{eq:dangerous_sufficient_dyadic}
(0.b_{n-1} b_n b_{n+1} \cdots)_2 \leq (0.01b_1 b_2 b_3 \cdots)_2 \enspace.
\end{equation}
\end{lemma}
Since $b_1 b_2 b_3 = 011$, the condition (\ref{eq:dangerous_sufficient_dyadic}) is satisfied if $b_{n-1} b_n = 00$, $b_{n-1} b_n b_{n+1} b_{n+2} = 0100$, or $b_{n-1} b_n b_{n+1} b_{n+2} b_{n+3} = 01010$.
Summarizing, we obtain the following sufficient condition for an accuracy parameter $n$ being undesirable:
\begin{proposition}
\label{prop:dangerous_sufficient}
In the above setting, a parameter $n \geq 2$ is undesirable if $b_{n-1} b_n = 00$, $b_{n-1} b_n b_{n+1} b_{n+2} = 0100$, or $b_{n-1} b_n b_{n+1} b_{n+2} b_{n+3} = 01010$.
\end{proposition}
\begin{remark}
\label{rem:not_necessary}
In general, the sufficient condition given by Proposition \ref{prop:dangerous_sufficient} is not necessary for a parameter $n \geq 2$ being undesirable.
More precisely, there exists a gap between the sufficient conditions in Proposition \ref{prop:dangerous_sufficient} and in Lemma \ref{lem:dangerous_sufficient_dyadic}.
For instance, $n = 65$ satisfies the condition (\ref{eq:dangerous_sufficient_dyadic}) but not the condition in Proposition \ref{prop:dangerous_sufficient}.
Thus a more precise study of the condition (\ref{eq:dangerous_sufficient_dyadic}) would provide a better result.
The author hopes that the condition (\ref{eq:dangerous_sufficient_dyadic}) can motivate further interesting arguments by its self-referential structure.
\end{remark}

\section{Occurrence Rates of the Three Patterns}
\label{sec:occurrence}

Motivated by Proposition \ref{prop:dangerous_sufficient}, in this section we investigate the asymptotic occurrence rate of the three patterns $00$, $0100$, and $01010$ in an infinite binary word.
The result will be used to estimate the asymptotic occurrence rate of undesirable parameters.

To formulate the problem, we introduce the following notations.
For a finite or infinite binary word $w = w_1 w_2 w_3 \cdots$ ($w_i \in \{0,1\}$), let $\ell(w)$ denote the length of $w$.
Let $P(w)$ denote the set of indices $i \geq 2$ in $w$ such that one of the following three conditions holds:
\begin{itemize}
\item $\ell(w) \geq i$ and $w_{i-1} w_i = 00$;
\item $\ell(w) \geq i+2$ and $w_{i-1} w_i w_{i+1} w_{i+2} = 0100$;
\item $\ell(w) \geq i+3$ and $w_{i-1} w_i w_{i+1} w_{i+2} w_{i+3} = 01010$.
\end{itemize}
By Proposition \ref{prop:dangerous_sufficient}, a parameter $n \geq 2$ is undesirable if $n \in P(b)$, where $b = b_1 b_2 b_3 \cdots$ is the fractional part of the dyadic expansion of $\sqrt{2}$ as an infinite binary word.
Let $w^{(k)}$ denote the initial subword of $w$ of length $k$.
Moreover, let $Z(w)$ denote the set of indices $i$ in $w$ such that $w_i = 0$.
Then our main theorem in this section shows relations between the quantities
\begin{equation}
\label{eq:defn_inf}
r_{\inf}(w) = \liminf_{n \to \infty} \frac{ |Z(w^{(n)})| }{ n }\enspace \mbox{ and }\enspace R_{\inf}(w) = \liminf_{n \to \infty} \frac{ |P(w^{(n)})| }{ n } \enspace,
\end{equation}
and relations between the quantities
\begin{equation}
\label{eq:defn_sup}
r_{\sup}(w) = \limsup_{n \to \infty} \frac{ |Z(w^{(n)})| }{ n }\enspace \mbox{ and }\enspace R_{\sup}(w) = \limsup_{n \to \infty} \frac{ |P(w^{(n)})| }{ n } \enspace.
\end{equation}
By using the above notations, we state the main theorem as follows:
\begin{theorem}
\label{thm:main}
For any infinite binary word $w = w_1 w_2 w_3 \cdots$, let $r_{\inf}(w)$, $r_{\sup}(w)$, $R_{\inf}(w)$, and $R_{\sup}(w)$ be defined in (\ref{eq:defn_inf}) and (\ref{eq:defn_sup}).
Then we have
\begin{equation}
\label{eq:thm_main}
\frac{ 5 r_{\inf}(w) - 2 }{ 3 } \leq R_{\inf}(w) \leq r_{\inf}(w) \enspace \mbox{ and } \enspace \frac{ 5 r_{\sup}(w) - 2 }{ 3 } \leq R_{\sup}(w) \leq r_{\sup}(w) \enspace.
\end{equation}
Moreover, for any real number $2/5 \leq r \leq 1$, there exists an infinite binary word $w$ such that $r_{\inf}(w) = r_{\sup}(w) = r$ and $R_{\inf}(w) = R_{\sup}(w) = (5r - 2)/3$, therefore the lower bounds are achieved.
Similarly, for any $0 \leq r \leq 1$, there exists an infinite binary word $w$ such that $r_{\inf}(w) = r_{\sup}(w) = r$ and $R_{\inf}(x) = R_{\sup}(x) = r$, therefore the upper bounds are achieved.
\end{theorem}
Note that the lower bounds of $R_{\inf}(w)$ and $R_{\sup}(w)$ become trivial if $r_{\inf}(w) < 2/5$ and $r_{\sup}(w) < 2/5$, respectively.
The proof of Theorem \ref{thm:main} will be given in Section \ref{sec:proof}.

Regarding the problem in Section \ref{sec:logistic}, by applying Theorem \ref{thm:main} to the above word $w = b$ of the fractional part of the dyadic expansion of $\sqrt{2}$, we have the following theorem:
\begin{theorem}
\label{thm:undesirable_bound}
In the above setting, let $d_N$ denote the number of the undesirable parameters $n \leq N$.
Then we have
\begin{equation}
\liminf_{N \to \infty} \frac{ d_N }{ N } \geq \frac{ 5 r_{\inf}(b) - 2 }{ 3 } \mbox{ and } \limsup_{N \to \infty} \frac{ d_N }{ N } \geq \frac{ 5 r_{\sup}(b) - 2 }{ 3 } \enspace.
\end{equation}
In particular, if $r_{\sup}(b) > 2/5$, then there exist infinitely many undesirable parameters.
\end{theorem}
As a result, the (lower bound of the) asymptotic occurrence rate of zeroes in the dyadic expansion of $\sqrt{2}$ yields a lower bound of the asymptotic occurrence rate of undesirable parameters.
Note that there has been the following long-standing conjecture:
\begin{conjecture}
\label{conj:simply_normal}
$\sqrt{2}$ is simply normal to base $2$; that is, the asymptotic occurrence rate of zeroes in the dyadic expansion of $\sqrt{2}$ is $1/2$ (i.e., $r_{\inf}(b) = r_{\sup}(b) = 1/2$ in the above notations).
\end{conjecture}
This conjecture reflects our naive intuition that the dyadic expansion of $\sqrt{2}$ looks very random.
There have been some further observations that sound positive for the conjecture.
For instance, Borel \cite{Bor09} proved that almost every real number (in terms of Lebesgue measure) is simply normal to base $2$ (more strongly, is normal to every base $q \geq 2$).
By combining Conjecture \ref{conj:simply_normal} to Theorem \ref{thm:undesirable_bound}, we obtain the following result, that is very likely to show a disadvantage of the PRNGs mentioned in Section \ref{sec:intro}:
\begin{corollary}
\label{cor:many_undesirable_parameters}
If Conjecture \ref{conj:simply_normal} is true, then the numbers $d_N$ of undesirable parameters $n \leq N$ satisfy that $\liminf_{N \to \infty} d_N / N \geq 1/6$.
\end{corollary}

\section{Proof of Main Theorem}
\label{sec:proof}

In this section, we give a proof of Theorem \ref{thm:main} in Section \ref{sec:occurrence}.
First, the upper bounds of $R_{\inf}(w)$ and $R_{\sup}(w)$ in (\ref{eq:thm_main}) follow simply from the fact that the map $i \mapsto i-1$ is an injection from $P(w)$ to $Z(w)$ for any finite binary word $w$.
For the remaining assertions, first in Section \ref{subsec:proof_lower_bound} we prove the lower bounds in (\ref{eq:thm_main}).
Secondly, in Section \ref{subsec:proof_lower_bound_best} we construct, for any $2/5 \leq r \leq 1$, an infinite binary word $w$ such that $r_{\inf}(w) = r_{\sup}(w) = r$ and $R_{\inf}(w) = R_{\sup}(w) = (5r - 2)/3$.
Finally, in Section \ref{subsec:proof_upper_bound_best} we construct, for any $0 \leq r \leq 1$, an infinite binary word $w$ such that $r_{\inf}(w) = r_{\sup}(w) = r$ and $R_{\inf}(w) = R_{\sup}(w) = r$.

\subsection{Proof of Lower Bounds}
\label{subsec:proof_lower_bound}

We prove that $R_{\inf}(w) \geq (5 r_{\inf}(w) - 2)/3$ and $R_{\sup}(w) \geq (5 r_{\sup}(w) - 2)/3$ for any infinite binary word $w$.
In the proof, we use the following notations.
For any (finite or infinite) word $w = w_1 w_2 w_3 \cdots$ and indices $1 \leq i \leq j \leq \ell(w)$, let $w_{\left[i,j\right]} = w_i w_{i+1} \cdots w_{j-1} w_j$.
Let $\emptyset$ denote the empty word.
Let $W_N$ denote the set of binary words of length $N$.
Let $\prec$ denote the lexicographic order on $W_n$ excluding equalities, for instance, we have $1011 \prec 1100$ and $0010 \not\prec 0010$.
For two words $w$ and $w'$, we write $w \subset w'$ if $w = w'_{\left[i,j\right]}$ for some indices $i \leq j$.
Let $w^j = w w \cdots w$ ($j$ repetition of $w$) for any integer $j \geq 0$.

The outline of our proof is as follows.
In the proof, we investigate the maximum value of the number $|Z(u)|$ of zeroes in $u \in W_N$ subject to the condition that $|P(u)|$ is upper bounded by a fixed value.
This will yield a relation between the quantities $|P(w^{(n)})|$ and $|Z(w^{(n)})|$ for each initial subword $w^{(n)}$ of a given infinite binary word $w$, from which the desired lower bounds will be derived.
For the purpose, we will introduce some \lq\lq elementary transformations'' for the words $u \in W_N$ that preserve $\ell(u)$ and $|Z(u)|$ and do not increase $|P(u)|$.
By iterating such elementary transformations, our argument will be reduced to the case of words in $W_N$ of some \lq\lq normal form'' that can be dealt with by case-by-case analysis.

We start the above program.
First, we introduce the following seven maps $\varphi_k:W_N \to W_N$, $1 \leq k \leq 7$, as aforementioned elementary transformations, where $v$ and $v'$ signify some (possibly empty) binary words.
We define
\begin{displaymath}
\varphi_1(u) = \begin{cases}
1^p v 0 & \mbox{, if } u = v 0 1^p,\ p \geq 1 \enspace; \\
u & \mbox{, otherwise}
\end{cases}
\end{displaymath}
(namely, $\varphi_1$ moves the ones at the tail of the word $u$, to the front of $u$; for instance, $\varphi_1(10100\underline{11}) = \underline{11}10100$ and $\varphi_1$ fixes $10100$),
\begin{displaymath}
\varphi_2(u) = \begin{cases}
1^{p+1} v 11 v' & \mbox{, if } u = 1^p v 111 v',\ p \geq 0,\ 111 \not\subset v \neq \emptyset,\ v_1 = v_{\ell(v)} = 0 \enspace; \\
u & \mbox{, otherwise}
\end{cases}
\end{displaymath}
(namely, $\varphi_2$ picks up a one from the first block of at least three consecutive ones after a zero and moves it to the front; for instance, $\varphi_2(\underline{1^6}0110\underline{1^4}01^50) = \underline{1^7}0110\underline{1^3}01^50$ and $\varphi_2$ fixes $1^30011010$),
\begin{displaymath}
\varphi_3(u) = \begin{cases}
v 0 11 0^{p-1} v' & \mbox{, if } u = v 0^p 11 v',\ p \geq 2,\ 0011 \not\subset v,\ v_{\ell(v)} \neq 0 \enspace; \\
u & \mbox{, otherwise}
\end{cases}
\end{displaymath}
(namely, $\varphi_3$ focuses on the first block of the form $0^p11$ with $p \geq 2$, and moves all but one zeroes in that block to the tail of that block; for instance, $\varphi_3(11011\underline{0^411}100110) = 11011\underline{0110^3}100110$ and $\varphi_3$ fixes $1011011$),
\begin{displaymath}
\varphi_4(u) = \begin{cases}
v 01100 v' & \mbox{, if } u = v 01010 v',\ 01010 \not\subset v010 \enspace; \\
u & \mbox{, otherwise}
\end{cases}
\end{displaymath}
(namely, $\varphi_4$ focuses on the first block of the form $01010$, and permutes the third and the fourth bits in that block; for instance, $\varphi_4(11\underline{01010}10) = 11\underline{01100}10$ and $\varphi_4$ fixes $011010110101$),
\begin{displaymath}
\varphi_5(u) = \begin{cases}
v 1 0^{p+2} v' & \mbox{, if } u = v 0^p 100 v',\ p \geq 1,\ 0100 \not\subset v 0^p,\ v_{\ell(v)} \neq 0 \enspace; \\
u & \mbox{, otherwise}
\end{cases}
\end{displaymath}
(namely, $\varphi_5$ focuses on the first block of the form $0^p100$ with $p \geq 1$, and moves the unique one in that block to the front of that block; for instance, $\varphi_5(10011\underline{0^3100}100) = 10011\underline{10^5}100$ and $\varphi_5$ fixes $100110010$),
\begin{displaymath}
\varphi_6(u) = \begin{cases}
v 01011 0^p v' & \mbox{, if } u = v 0^p 10110 v',\ p \geq 2,\ 0010110 \not\subset v 0^p,\ v_{\ell(v)} \neq 0; \\
u & \mbox{, otherwise}
\end{cases}
\end{displaymath}
(namely, $\varphi_6$ focuses on the first block of the form $0^p10110$ with $p \geq 2$, and moves all but one zeroes at the beginning of that block to the tail of that block; for instance, $\varphi_6(1\underline{0^410110}0101100) = 1\underline{010110^4}0101100$ and $\varphi_6$ fixes $1010110$), and
\begin{displaymath}
\varphi_7(u) = \begin{cases}
v 1010110 v' & \mbox{, if } u = v 0110110 v',\ 0110110 \not\subset v 0110 \enspace; \\
u & \mbox{, otherwise}
\end{cases}
\end{displaymath}
(namely, $\varphi_7$ focuses on the first block of the form $0110110$, and permutes the first and the second bits in that block; for instance, $\varphi_7(1^30\underline{0110110}110) = 1^30\underline{1010110}110$ and $\varphi_7$ fixes $0111011010$).

Let $W_N^{\varphi}$ denote the set of all $u \in W_N$ that are fixed by every $\varphi_k$, $1 \leq k \leq 7$.
Note that each of the seven maps $\varphi_k$ is well-defined and satisfies that $\ell(\varphi_k(u)) = \ell(u)$ and $|Z(\varphi_k(u))| = |Z(u)|$, since $\varphi_k$ is just a permutation of bits in a given word.
Moreover, it follows immediately from the definition that each $\varphi_k$ is a weakly increasing map with respect to $\prec$, namely we have $u \preceq \varphi_k(u)$.
Since $W_N$ is a finite set, this implies that each $u \in W_N$ can be transformed to a word $\overline{u} \in W_N^{\varphi}$ by finitely many times of applications of the maps $\varphi_k$, $1 \leq k \leq 7$.
Note that this $\overline{u}$ is not necessarily unique for a given $u \in W_N$ due to various choices of the order of applying the maps $\varphi_k$.

To reduce our argument to the case of the words in $W_N^{\varphi}$, we would like to show that $|P(\overline{u})| \leq |P(u)|$ for any $u \in W_N$.
For the purpose, it suffices to show that $|P(\varphi_k(u))| \leq |P(u)|$ for every map $\varphi_k$, $1 \leq k \leq 7$.
This is proceeded by the following seven lemmas, where we use the notation:
\begin{displaymath}
P_{i,j}(u) = P(u) \cap \{i,i+1,\dots,j-1,j\} \mbox{ for any indices } i \leq j \mbox{ in } u \in W_N \enspace.
\end{displaymath}
Before giving the lemmas, note that for any word $u$ and any index $i$,
\begin{equation}
\label{eq:property_1}
\mbox{we have } i \not\in P(u) \mbox{ unless } i \geq 2 \mbox{ and } u_{i-1} = 0 \enspace,
\end{equation}
therefore $1 \not\in P(u)$.
Similarly,
\begin{equation}
\label{eq:property_2}
\mbox{if } u_i = 1, \mbox{ then we have } i \not\in P(u) \mbox{ unless } 2 \leq i \leq \ell(u) - 1 \mbox{ and } u_{i-1} = u_{i+1} = 0 \enspace.
\end{equation}
Moreover, it is obvious that
\begin{equation}
\label{eq:property_3}
\mbox{whether } i \in P(u) \mbox{ or not depends solely on } u_{\left[i-1,i+3\right]} \enspace.
\end{equation}
Now we show the lemmas as follows, where we write $u' = \varphi_k(u)$ for the map $\varphi_k$ under consideration:
\begin{lemma}
\label{lem:transform_1}
If $u \in W_N$, then $|P(\varphi_1(u))| = |P(u)|$.
\end{lemma}
\begin{proof}
It suffices to consider the case that $u' \neq u$, namely $u = v 0 1^p$ and $u' = 1^p v 0$ with $p \geq 1$, as in the former case of the definition of $\varphi_1$.
Now if $x = u'_{\left[i-1,j\right]}$ is a subword in $u'$ of one of the three forms $00$, $0100$, or $01010$, corresponding to an index $i \in P(u')$, then $x$ should be contained in $v0$ by the shapes of $x$ and $u'$, therefore $u_{\left[i-1-p,j-p\right]} = x$ and $i-p \in P(u)$.
Similarly, if $x = u_{\left[i-1,j\right]}$ is a subword in $u$ of the form $00$, $0100$, or $01010$, corresponding to an $i \in P(u)$, then $x \subset v0$, therefore $u'_{\left[i-1+p,j+p\right]} = x$ and $i+p \in P(u')$.
Thus $i \mapsto i + p$ is a bijection $P(u) \to P(u')$, therefore Lemma \ref{lem:transform_1} holds.
\end{proof}
\begin{lemma}
\label{lem:transform_2}
If $u \in W_N$, then $|P(\varphi_2(u))| = |P(u)|$.
\end{lemma}
\begin{proof}
It suffices to consider the case that $u = 1^p v 111 v'$ and $u' = 1^{p+1} v 11 v'$ as in the former case of the definition of $\varphi_2$.
Now by the shapes of $u$ and $u'$, any subword in $u$ of the form $00$, $0100$, or $01010$ is contained in either $v$ or $v'$, and the same also holds for $u'$.
This implies that there exists a bijection $P(u) \to P(u')$, hence Lemma \ref{lem:transform_2} holds.
\end{proof}
\begin{lemma}
\label{lem:transform_3}
If $u \in W_N$, then $|P(\varphi_3(u))| \leq |P(u)|$.
\end{lemma}
\begin{proof}
It suffices to consider the case that $u = v 0^p 11 v'$ and $u' = v 0110^{p-1} v'$ as in the former case of the definition of $\varphi_3$.
Put $\ell = \ell(v)$.
Then for any subword $x$ in $u'$ of the form $00$, $0100$, or $01010$ corresponding to an $i \in P(u')$, one of the following four conditions is satisfied:
\begin{enumerate}
\item $x$ is contained in the block $v0$;
\item $x = 00$ and $x$ is contained in the block $0^{p-1}$ (thus $\ell + 5 \leq i \leq \ell + p + 2$);
\item $i = \ell + p + 3$, namely $i$ is the first position of the block $v'$;
\item $x$ is contained in the block $v'$.
\end{enumerate}
In the cases 1 and 4, $x$ is also contained in $u$ and we have $i \in P(u)$.
In the case 2, $x$ is also contained in the last $p-1$ bits of the block $0^p$ in $u$, and we have $i-2 \in P(u)$.
Moreover, we have $\ell + 2 \in P(u)$ since $p \geq 2$.
Thus there exists an injection $P(u') \to P(u)$ that maps $i \in P(u')$ to $i$ for the cases 1 and 4, to $i-2$ for the case 2, and to $\ell + 2$ for the case 3.
Hence Lemma \ref{lem:transform_3} holds.
\end{proof}
\begin{lemma}
\label{lem:transform_4}
If $u \in W_N$, then $|P(\varphi_4(u))| \leq |P(u)|$.
\end{lemma}
\begin{proof}
It suffices to consider the case that $u = v 01010 v'$ and $u' = v 01100 v'$ as in the former case of the definition of $\varphi_4$.
Put $\ell = \ell(v)$.
Then for any subword $x$ in $u'$ of the form $00$, $0100$, or $01010$ corresponding to an $i \in P(u')$, one of the following three conditions is satisfied:
\begin{enumerate}
\item $x$ is contained in the block $v0$;
\item $i = \ell + 5$, namely $i$ is the last position of the block $01100$;
\item $x$ is contained in the block $0v'$ (thus $\ell + 6 \leq i$).
\end{enumerate}
In the cases 1 and 3, $x$ is also contained in $u$ and we have $i \in P(u)$.
Since $\ell + 2 \in P(u)$, there exists an injection $P(u') \to P(u)$ that maps $i \in P(u')$ to $i$ for the cases 1 and 3, and to $\ell + 2$ for the case 2.
Hence Lemma \ref{lem:transform_4} holds.
\end{proof}
\begin{lemma}
\label{lem:transform_5}
If $u \in W_N$, then $|P(\varphi_5(u))| \leq |P(u)|$.
\end{lemma}
\begin{proof}
It suffices to consider the case that $u = v 0^p 100 v'$ and $u' = v 10^{p+2} v'$ as in the former case of the definition of $\varphi_5$.
Put $\ell = \ell(v)$.
Note that $p \geq 1$ and $v_{\ell} \neq 0$ by the definition of $\varphi_5$.
Then for any subword $x$ in $u'$ of the form $00$, $0100$, or $01010$ corresponding to an $i \in P(u')$, one of the following four conditions is satisfied:
\begin{enumerate}
\item $x$ is contained in the block $v$;
\item $x$ is contained in the first $p$ bits of the block $0^{p+2}$ (thus $\ell + 3 \leq i \leq \ell + p + 1$);
\item $i = \ell + p + 2$, namely $i$ is the second last position of the block $0^{p+2}$;
\item $x$ is contained in the block $00v'$.
\end{enumerate}
In the cases 1 and 4, $x$ is also contained in $u$ and we have $i \in P(u)$.
In the case 2, $x$ is also contained in the block $0^p$ in $u$, and we have $i-1 \in P(u)$.
Moreover, we have $\ell + p + 1 \in P(u)$ since $p \geq 1$ (namely $u_{\left[\ell + p,\ell + p + 3\right]} = 0100$).
Thus there exists an injection $P(u') \to P(u)$ that maps $i \in P(u')$ to $i$ for the cases 1 and 4, to $i-1$ for the case 2, and to $\ell + p + 1$ for the case 3.
Hence Lemma \ref{lem:transform_5} holds.
\end{proof}
\begin{lemma}
\label{lem:transform_6}
If $u \in W_N$, then $|P(\varphi_6(u))| \leq |P(u)|$.
\end{lemma}
\begin{proof}
It suffices to consider the case that $u = v 0^p 10110 v'$ and $u' = v 010110^p v'$ as in the former case of the definition of $\varphi_6$.
Put $\ell = \ell(v)$.
Note that $p \geq 2$ and $v_{\ell} \neq 0$ by the definition of $\varphi_6$.
Then for any subword $x$ in $u'$ of the form $00$, $0100$, or $01010$ corresponding to an $i \in P(u')$, one of the following four conditions is satisfied:
\begin{enumerate}
\item $x$ is contained in the block $v0$ (thus $i \leq \ell - 1$ since $v_{\ell} \neq 0$);
\item $\ell \geq 2$, $v_{\left[\ell-1,\ell\right]} = 01$, $x = 01010$ and $i = \ell$;
\item $x$ is contained in the block $0^p$ (thus $\ell + 7 \leq i \leq \ell + p + 5$);
\item $x$ is contained in the block $0v'$ (thus $\ell + p + 6 \leq i$).
\end{enumerate}
In the cases 1 and 4, $x$ is also contained in $u$ and we have $i \in P(u)$.
In the case 3, $x$ is also contained in the block $0^p$ in $u$, and we have $i-5 \in P(u)$.
Moreover, in the case 2, we have $\ell \in P(u)$ since $p \geq 2$ (namely $u_{\left[\ell-1,\ell+2\right]} = 0100$).
Thus there exists an injection $P(u') \to P(u)$ that maps $i \in P(u')$ to $i$ for the cases 1 and 4, to $i-5$ for the case 3, and to $\ell$ for the case 2.
Hence Lemma \ref{lem:transform_6} holds.
\end{proof}
\begin{lemma}
\label{lem:transform_7}
If $u \in W_N$, then $|P(\varphi_7(u))| \leq |P(u)|$.
\end{lemma}
\begin{proof}
It suffices to consider the case that $u = v 0110110 v'$ and $u' = v 1010110 v'$ as in the former case of the definition of $\varphi_7$.
Put $\ell = \ell(v)$.
Then for any subword $x$ in $u'$ of the form $00$, $0100$, or $01010$ corresponding to an $i \in P(u')$, one of the following four conditions is satisfied:
\begin{enumerate}
\item $x$ is contained in the block $v$;
\item $\ell \geq 1$, $v_{\ell} = 0$, $x = 01010$ and $i = \ell + 1$;
\item $x$ is contained in the block $0v'$ (thus $\ell + 8 \leq i$).
\end{enumerate}
In the cases 1 and 3, $x$ is also contained in $u$ and we have $i \in P(u)$.
Moreover, in the case 2, we have $\ell + 1 \in P(u)$ (namely $u_{\left[\ell,\ell+1\right]} = 00$).
Thus there exists an injection $P(u') \to P(u)$ that maps $i \in P(u')$ to $i$ for the cases 1 and 3, and to $\ell + 1$ for the case 2.
Hence Lemma \ref{lem:transform_7} holds.
\end{proof}
Thus we have proven that $|P(\overline{u})| \leq |P(u)|$ for any $u \in W_N$ as desired.
From now, we determine the possibilities of the shape of $\overline{u} \in W_N^{\varphi}$.
For the purpose, first we show that any word in $W_N^{\varphi}$ does not contain a subword of type 1--11 in Table \ref{tab:excluded_patterns}.
For instance, if $\overline{u} \in W_N^{\varphi}$, then we have $010111 \not\subset \overline{u}$ since $010111$ is a word of type 2.
In fact, only subwords of types 1, 2, 4, 5, 7, and 11 will appear in the subsequent argument.
However, we also include other subwords in Table \ref{tab:excluded_patterns} since these are used in the proof of the above fact (Lemma \ref{lem:excluded_patterns} below).
Now we show the following lemma:
\begin{table}[htb]
\centering
\caption{Excluded subwords for words in $W_N^{\varphi}$}
\label{tab:excluded_patterns}
Here $v$ is a (possibly empty) word, and \lq $)$' for type 1 denotes the tail of the word $\overline{u}$.\\
\begin{tabular}{|c|c||c|c||c|c|} \hline
type & subword & type & subword & type & subword \\ \hline\hline
1 & $0 v 1 )$ & 2 & $0 v 111$ & 3 & $0011$ \\ \hline
4 & $01010$ & 5 & $0100$ & 6 & $0010110$ \\ \hline
7 & $0110110$ & 8 & $001011$ & 9 & $00101$ \\ \hline
10 & $0010v$ ($v \neq \emptyset$) & 11 & $001 v$ ($v \neq 0$) & \multicolumn{2}{|c}{} \\ \cline{1-4}
\end{tabular}
\end{table}
\begin{lemma}
\label{lem:excluded_patterns}
Any $\overline{u} \in W_N^{\varphi}$ does not contain a subword listed in Table \ref{tab:excluded_patterns}.
\end{lemma}
\begin{proof}
First, we show that for each $1 \leq k \leq 7$, we have $\varphi_k(w) \neq w$ if $w \in W_N$ contains a subword $x$ of type $k$ in Table \ref{tab:excluded_patterns}.
It suffices to show that $w$ satisfies the condition in the first case of the definition of $\varphi_k$.
In the case $k = 1$, by focusing on the last zero in $x$ and the subsequent part of $x$ (the former exists and the latter is nonempty by the shape of $x$), it follows that $w$ ends with $01^p$, $p \geq 1$, therefore the claim holds.
In the case $k = 2$, the shape of $x$ implies that $w$ contains a subword of the form $0111$ (consider the last zero in $x$ and the subsequent part of $x$), and by focusing on the leftmost such subword in $w$ the claim follows.
In the case $k = 3$, the shape of $x$ implies that $w$ contains a subword of the form $0^p 11$, $p \geq 2$, and by focusing on the leftmost such subword in $w$ the claim follows.
In the case $k = 4$, the claim follows by focusing on the leftmost subword in $w$ of the form $01010$ (the same as $x$).
In the case $k = 5$, the shape of $x$ implies that $w$ contains a subword of the form $0^p 100$, $p \geq 1$, and by focusing on the leftmost such subword in $w$ the claim follows.
In the case $k = 6$, the shape of $x$ implies that $w$ contains a subword of the form $0^p 10110$, $p \geq 2$, and by focusing on the leftmost such subword in $w$ the claim follows.
Finally, in the case $k = 7$, the claim follows by focusing on the leftmost subword in $w$ of the form $0110110$ (the same as $x$).
Thus we have shown that any $\overline{u} \in W_N^{\varphi}$ does not contain a subword of types 1--7 in Table \ref{tab:excluded_patterns}.

For subwords of type 8, if $001011 \subset \overline{u}$, then $\overline{u}$ must contain one of the three subwords $0010110$, $0010111$, and $001011)$ that are of types 6, 2, and 1 in Table \ref{tab:excluded_patterns}, respectively.
This contradicts the previous paragraph, therefore $001011 \not\subset \overline{u}$ as desired.

For subwords of type 9, if $00101 \subset \overline{u}$, then $\overline{u}$ must contain one of the three subwords $0\underline{01010}$ (that contains a subword of type 4 in Table \ref{tab:excluded_patterns}), $001011$ (a subword of type 8 in Table \ref{tab:excluded_patterns}), and $00101)$ (a subword of type 1 in Table \ref{tab:excluded_patterns}).
This contradicts the previous paragraphs, therefore $00101 \not\subset \overline{u}$ as desired.

For subwords of type 10, if $0010v \subset \overline{u}$ for some word $v \neq \emptyset$, then $\overline{u}$ must contain one of the two subwords $0\underline{0100}$ (that contains a subword of type 5 in Table \ref{tab:excluded_patterns}) and $00101$ (a subword of type 9 in Table \ref{tab:excluded_patterns}).
This contradicts the previous paragraphs, therefore $0010v \not\subset \overline{u}$ as desired.

Finally, for subwords of type 11, if $001v \subset \overline{u}$ for some word $v \neq 0$, then $\overline{u}$ must contain one of the three subwords $0010v'$ with $v' \neq \emptyset$ (a subword of type 10 in Table \ref{tab:excluded_patterns}), $0011$ (a subword of type 3 in Table \ref{tab:excluded_patterns}), or $001)$ (a subword of type 1 in Table \ref{tab:excluded_patterns}).
This contradicts the previous paragraphs, therefore $001v \not\subset \overline{u}$ as desired.
Hence Lemma \ref{lem:excluded_patterns} holds.
\end{proof}
Owing to Lemma \ref{lem:excluded_patterns}, we obtain the following classification of the words in $W_N^{\varphi}$:
\begin{lemma}
\label{lem:expression_fixed_point}
Any word $\overline{u}$ in $W_N^{\varphi}$ is of one of the seven types in Table \ref{tab:complete_list_fixed_point}.
\end{lemma}
\begin{table}[htb]
\centering
\caption{Classification of words $u$ in $W_N^{\varphi}$}
\label{tab:complete_list_fixed_point}
\begin{tabular}{|c|cc|} \hline
Type 1 & \multicolumn{2}{|c|}{$u = 1^p 0^q$\quad ($p \geq 0$, $q \geq 0$)} \\ \cline{2-3}
& \multicolumn{2}{|c|}{$N = p + q$} \\
& $|Z(u)| = q$ & $|P(u)| = q - 1$ \\ \cline{2-3}
& \multicolumn{2}{|c|}{$|Z(u)|/N = |P(u)| / N + 1/N$} \\ \hline\hline
Type 2 & \multicolumn{2}{|c|}{$u = 1^p (01011)^s 0^q 10$\quad ($p \geq 0$, $q \geq 2$, $s \geq 0$)} \\ \cline{2-3}
& \multicolumn{2}{|c|}{$N = 5s + p + q + 2$} \\
& $|Z(u)| = 2s + q + 1$ & $|P(u)| = q - 1$ \\ \cline{2-3}
& \multicolumn{2}{|c|}{$|Z(u)| / N = 3/5 \cdot |P(u)| / N + 2/5 + 4/(5N) - 2p/(5N)$} \\ \hline\hline
Type 3 & \multicolumn{2}{|c|}{$u = 1^p 011(01011)^s 0^q 10$\quad ($p \geq 0$, $q \geq 2$, $s \geq 0$)} \\ \cline{2-3}
& \multicolumn{2}{|c|}{$N = 5s + p + q + 5$} \\
& $|Z(u)| = 2s + q + 2$ & $|P(u)| = q - 1$ \\ \cline{2-3}
& \multicolumn{2}{|c|}{$|Z(u)| / N = 3/5 \cdot |P(u)| / N + 2/5 + 3/(5N) - 2p/(5N)$} \\ \hline\hline
Type 4 & \multicolumn{2}{|c|}{$u = 1^p (01011)^s 0^q$\quad ($p \geq 0$, $q \geq 1$, $s \geq 1$)} \\ \cline{2-3}
& \multicolumn{2}{|c|}{$N = 5s + p + q$} \\
& $|Z(u)| = 2s + q$ & $|P(u)| = q - 1$ \\ \cline{2-3}
& \multicolumn{2}{|c|}{$|Z(u)| / N = 3/5 \cdot |P(u)| / N + 2/5 + 3/(5N) - 2p/(5N)$} \\ \hline\hline
Type 5 & \multicolumn{2}{|c|}{$u = 1^p 011(01011)^s 0^q$\quad ($p \geq 0$, $q \geq 1$, $s \geq 0$)} \\ \cline{2-3}
& \multicolumn{2}{|c|}{$N = 5s + p + q + 3$} \\
& $|Z(u)| = 2s + q + 1$ & $|P(u)| = q - 1$ \\ \cline{2-3}
& \multicolumn{2}{|c|}{$|Z(u)| / N = 3/5 \cdot |P(u)| / N + 2/5 + 2/(5N) - 2p/(5N)$} \\ \hline\hline
Type 6 & \multicolumn{2}{|c|}{$u = 1^p (01011)^s 010$\quad ($p \geq 0$, $s \geq 0$)} \\ \cline{2-3}
& \multicolumn{2}{|c|}{$N = 5s + p + 3$} \\
& $|Z(u)| = 2s + 2$ & $|P(u)| = 0$ \\ \cline{2-3}
& \multicolumn{2}{|c|}{$|Z(u)| / N = 2/5 + 4/(5N) - 2p/(5N)$} \\ \hline\hline
Type 7 & \multicolumn{2}{|c|}{$u = 1^p 011(01011)^s 010$\quad ($p \geq 0$, $s \geq 0$)} \\ \cline{2-3}
& \multicolumn{2}{|c|}{$N = 5s + p + 6$} \\
& $|Z(u)| = 2s + 3$ & $|P(u)| = 0$ \\ \cline{2-3}
& \multicolumn{2}{|c|}{$|Z(u)| / N = 2/5 + 3/(5N) - 2p/(5N)$} \\ \hline
\end{tabular}
\end{table}
\begin{proof}
First, note that any $\overline{u} \in W_N^{\varphi}$ can be expressed in the following form:
\begin{displaymath}
\overline{u} = 1^{p_0} 0^{q_1} 1^{p_1} \cdots 0^{q_k} 1^{p_k} 0^{q_{k+1}},\ 
k \geq 0,\ p_0 \geq 0,\ q_{k+1} \geq 0,\ p_i \geq 1,\ q_i \geq 1 \mbox{ ($1 \leq i \leq k$)} \enspace.
\end{displaymath}
We apply Lemma \ref{lem:excluded_patterns} to this $\overline{u}$.
First, the absence of a subword of type 1 in Table \ref{tab:excluded_patterns} implies that $\overline{u}$ does not end with a one unless $\overline{u}$ contains no zeroes.
Thus we have $q_{k+1} \geq 1$ if $k \geq 1$.
Secondly, the absence of a subword of type 2 in Table \ref{tab:excluded_patterns} implies that three consecutive ones do not appear after a zero, therefore we have $p_i \in \{1,2\}$ for every $1 \leq i \leq k$.
Moreover, the absence of a subword of type 11 in Table \ref{tab:excluded_patterns} implies that if $001 \subset \overline{u}$, then a zero should follow that subword $001$ immediately and $\overline{u}$ should end with that zero.
By these conditions, the possible shapes of $\overline{u}$ are classified as follows:
\begin{enumerate}
\item $\overline{u} = 1^{p_0} 0^{q_1}$ (corresponding to the case $k = 0$);
\item $\overline{u} = 1^{p_0} 0 1^{p_1} 0^{q_2}$, $p_1 \in \{1,2\}$, $q_2 \geq 1$ (corresponding to the case $k = 1$, $q_k = 1$);
\item $\overline{u} = 1^{p_0} 0^{q_1} 10$, $q_1 \geq 2$ (corresponding to the case $k = 1$, $q_k \geq 2$);
\item $\overline{u} = 1^{p_0} 0 1^{p_1} \cdots 0 1^{p_{k-1}} 0 1^{p_k} 0^{q_{k+1}}$, $p_i \in \{1,2\}$ ($1 \leq i \leq k$), $q_{k+1} \geq 1$ (corresponding to the case $k \geq 2$, $q_k = 1$);
\item $\overline{u} = 1^{p_0} 0 1^{p_1} \cdots 0 1^{p_{k-1}} 0^{q_k} 10$, $p_i \in \{1,2\}$ ($1 \leq i \leq k$), $q_k \geq 2$ (corresponding to the case $k \geq 2$, $q_k \geq 2$).
\end{enumerate}
Case 1 corresponds to Type 1 in Table \ref{tab:complete_list_fixed_point}.
In Case 2, a choice $p_1 = 1$ implies that $q_2 = 1$ by the absence of a subword $0100$ of type 5 in Table \ref{tab:excluded_patterns}, and it corresponds to Type 6 in Table \ref{tab:complete_list_fixed_point} with parameter $s = 0$.
On the other hand, the other choice $p_1 = 2$ corresponds to Type 5 in Table \ref{tab:complete_list_fixed_point} with parameter $s = 0$.
Case 3 corresponds to Type 2 in Table \ref{tab:complete_list_fixed_point} with parameter $s = 0$.

The remaining part of the proof focuses on Cases 4 and 5.
In Case 4, the absence of a subword $01010$ of type 4 in Table \ref{tab:excluded_patterns} and a subword $0110110$ of type 7 in Table \ref{tab:excluded_patterns} implies that $(p_i,p_{i+1}) = (1,2)$ or $(2,1)$ for each $1 \leq i \leq k-1$.
Thus the sequence $(p_1,p_2,\dots,p_k)$ is of one of the four forms $(1,2,1,2,\dots,1,2)$, $(1,2,1,2,\dots,2,1)$, $(2,1,2,1,\dots,2,1)$, and $(2,1,2,1,\dots,1,2)$.
The first and the fourth cases correspond to Type 4 and Type 5 in Table \ref{tab:complete_list_fixed_point}, respectively.
On the other hand, the second and the third cases correspond to Type 6 and Type 7 in Table \ref{tab:complete_list_fixed_point}, respectively, since now we have $p_k = 1$ and the absence of a subword $0100$ of type 5 in Table \ref{tab:excluded_patterns} implies that $q_{k+1} = 1$.

Finally, in Case 5, the fact $q_k \geq 2$ and the absence of a subword $0100$ of type 5 in Table \ref{tab:excluded_patterns} imply that $p_{k-1} = 2$.
Now by the same argument as the previous paragraph, the sequence $(p_1,p_2,\dots,p_{k-1})$ is either $(1,2,1,2,\dots,1,2)$, or $(2,1,2,1,\dots,1,2)$.
These cases correspond to Type 2 and Type 3 in Table \ref{tab:complete_list_fixed_point}, respectively.
Hence Lemma \ref{lem:expression_fixed_point} holds.
\end{proof}
Table \ref{tab:complete_list_fixed_point} also includes, for each $u \in W_N^{\varphi}$, relations for the values $N$, $|Z(u)|$, $|P(u)|$, $p$, $q$ (except for Types 6 and 7), and $s$ (except for Type 1).

From now, we present the main part of the proof of lower bounds in Theorem \ref{thm:main}.
First we show that $R_{\inf}(w) \geq (5 r_{\inf}(w) - 2) / 3$ for any infinite binary word $w$.
This is trivial if $r_{\inf}(w) \leq 2/5$, therefore we focus on the case $r_{\inf}(w) > 2/5$.
Now we associate a (not necessarily unique) word $y(n) = \overline{w^{(n)}}$ in $W_n^{\varphi}$ to each initial subword $w^{(n)} \in W_n$ of $w$ by applying the maps $\varphi_i$, $1 \leq i \leq 7$, repeatedly.
Note that $|Z(y(n))| = |Z(w^{(n)})|$ and $|P(y(n))| \leq |P(w^{(n)})|$ by the above argument.
Now we have the following:
\begin{lemma}
\label{lem:proof_theorem_1}
$y(n) \in W_n^{\varphi}$ is not of type 6 or 7 in Table \ref{tab:complete_list_fixed_point} for any sufficiently large $n$.
\end{lemma}
\begin{proof}
Assume contrary that $y(n)$ is of type 6 or 7 for infinitely many $n$.
Then we have $|Z(w^{(n)})|/n = |Z(y(n))|/n \leq 2/5 + 4/(5n)$ for those $n$ by Table \ref{tab:complete_list_fixed_point}.
Note that the right-hand side converges to $2/5$ when $n \to \infty$.
On the other hand, by the assumption $r_{\inf}(w) > 2/5$, there exists an integer $N$ and a constant $c > 0$ such that $|Z(w^{(n)})| / n > 2/5 + c$ for any $n \geq N$.
This is a contradiction.
Hence Lemma \ref{lem:proof_theorem_1} holds.
\end{proof}
If $y(n)$ is of type 1 in Table \ref{tab:complete_list_fixed_point}, then we have
\begin{equation}
\label{eq:proof_theorem_type1}
\frac{ |P(w^{(n)})| }{ n } \geq \frac{ |P(y(n))| }{ n }
= \frac{ |Z(y(n))| }{ n } - \frac{ 1 }{ n }
= \frac{ |Z(w^{(n)})| }{ n } - \frac{ 1 }{ n } \enspace.
\end{equation}
On the other hand, if $y(n)$ is of types 2--5 in Table \ref{tab:complete_list_fixed_point}, then we have
\begin{equation}
\label{eq:proof_theorem_type2-5}
\begin{split}
\frac{ |P(w^{(n)})| }{ n } \geq \frac{ |P(y(n))| }{ n }
&= \frac{ 5 }{ 3 } \cdot \frac{ |Z(y(n))| }{ n } - \frac{ 2 }{ 3 } - \frac{ c }{ 3n } + \frac{ 2p }{ 3n } \\
&\geq \frac{ 5 }{ 3 } \cdot \frac{ |Z(w^{(n)})| }{ n } - \frac{ 2 }{ 3 } - \frac{ 4 }{ 3n } \enspace,
\end{split}
\end{equation}
where $c = 4$, $3$, $3$, and $2$ in the case of types 2, 3, 4, and 5, respectively.
Now we have $|Z(w^{(n)})| / n \leq 1$ by the definition, therefore the right-hand sides of (\ref{eq:proof_theorem_type1}) and (\ref{eq:proof_theorem_type2-5}) are larger than or equal to $5/3 \cdot |Z(w^{(n)})| / n - 2/3 - 4/(3n)$.
This implies that
\begin{equation}
\label{eq:proof_theorem_alltype}
\frac{ |P(w^{(n)})| }{ n } \geq \frac{ 5 }{ 3 } \cdot \frac{ |Z(w^{(n)})| }{ n } - \frac{ 2 }{ 3 } - \frac{ 4 }{ 3n }
\end{equation}
for any $n$ such that $y(n)$ is of types 1--5.
By Lemma \ref{lem:proof_theorem_1}, any sufficiently large $n$ satisfies the condition.
Now the desired bound $R_{\inf}(w) \geq (5 r_{\inf}(w) - 2) / 3$ is derived by taking the $\liminf_{n \to \infty}$ of both sides of (\ref{eq:proof_theorem_alltype}).

Secondly, we show that $R_{\sup}(w) \geq (5 r_{\sup}(w) - 2) / 3$ for any infinite binary word $w$.
This is trivial if $r_{\sup}(w) \leq 2/5$, therefore we focus on the case that $r_{\sup}(w) > 2/5$.
We define $y(n) \in W_n^{\varphi}$ for each $n$ in the same way as above.
Now it suffices to show that, for any $\varepsilon > 0$, there exist infinitely many indices $n$ such that $|P(w^{(n)})| / n > (5 r_{\sup}(w) - 2)/3 - \varepsilon$.
Take an $\varepsilon'$ such that $0 < \varepsilon' < 3 \varepsilon / 10$ and $\varepsilon' < r_{\sup}(w) - 2/5$.
Then by the definition of $r_{\sup}(w)$, there exist infinitely many indices $n$ such that $|Z(w^{(n)})| / n > r_{\sup}(w) - \varepsilon'$; let $\mathcal{N}$ denote the (infinite) set of such indices $n$.
Now we have the following:
\begin{lemma}
\label{lem:proof_theorem_2}
$y(n) \in W_n^{\varphi}$ is not of type 6 or 7 in Table \ref{tab:complete_list_fixed_point} for any sufficiently large $n \in \mathcal{N}$.
\end{lemma}
\begin{proof}
If $n \in \mathcal{N}$ and $y(n)$ is of type 6 or 7, then we have $|Z(y(n))| / n \leq 2/5 + 4/(5n)$ by Table \ref{tab:complete_list_fixed_point}, while $|Z(y(n))| / n = |Z(w^{(n)})| / n > r_{\sup}(w) - \varepsilon'$ by the definition of $\mathcal{N}$.
Thus we have $4/(5n) > r_{\sup}(w) - \varepsilon' - 2/5$ for such an $n$.
On the other hand, we have $r_{\sup}(w) - \varepsilon' - 2/5 > 0$ by the choice of $\varepsilon'$.
This implies that the condition for $n$ is not satisfied when $n \in \mathcal{N}$ is sufficiently large.
Hence Lemma \ref{lem:proof_theorem_2} holds.
\end{proof}
By Lemma \ref{lem:proof_theorem_2}, we may assume without loss of generality that $y(n)$ is of Types 1--5 for any $n \in \mathcal{N}$.

The relation (\ref{eq:proof_theorem_alltype}) above also holds in this case if $y(n)$ is of types 1--5 in Table \ref{tab:complete_list_fixed_point}.
By Lemma \ref{lem:proof_theorem_2}, any sufficiently large $n \in \mathcal{N}$ satisfies this condition.
For those (infinitely many) $n$, the definition of $\mathcal{N}$ and the choice of $\varepsilon'$ imply that
\begin{equation}
\label{eq:proof_theorem_alltype_2}
\frac{ |P(w^{(n)})| }{ n } > \frac{ 5 }{ 3 } (r_{\sup}(w) - \varepsilon') - \frac{ 2 }{ 3 } - \frac{ 4 }{ 3n }
> \frac{ 5r_{\sup}(w) - 2 }{ 3 } - \frac{ \varepsilon }{ 2 } - \frac{ 4 }{ 3n } \enspace.
\end{equation}
Moreover, the right-hand side of (\ref{eq:proof_theorem_alltype_2}) is larger than $(5r_{\sup}(w) - 2)/3 - \varepsilon$ if $n$ is sufficiently large.
Thus we have $|P(w^{(n)})| / n > (5r_{\sup}(w) - 2)/3 - \varepsilon$ for infinitely many $n$.
Since $\varepsilon > 0$ is arbitrary, it follows that $R_{\sup}(w) \geq (5 r_{\sup}(w) - 2) / 3$ as desired.
Hence the proof of lower bounds in Theorem \ref{thm:main} is concluded.

\subsection{Construction of Words Achieving Lower Bounds}
\label{subsec:proof_lower_bound_best}

We show that for any $2/5 \leq r \leq 1$, there exists an infinite binary word $w$ such that $r_{\inf}(w) = r_{\sup}(w) = r$ and $R_{\inf}(w) = R_{\sup}(w) = (5r - 2) / 3$.
First, if $r = 2/5$, then $w = 0101101011 \cdots$ (infinite repetition of $01011$) satisfies the condition.
On the other hand, if $r = 1$, then $w = 0000 \cdots$ (infinite repetition of $0$) satisfies the condition.
From now, we focus on the remaining case $2/5 < r < 1$.

We introduce some auxiliary notations.
First, put
\begin{displaymath}
p = \left\lceil \frac{ 5r-2 }{ 1-r } \right\rceil \mbox{ and } \alpha = p + 5 - \frac{ 3 }{ 1-r } = p - \frac{ 5r-2 }{ 1-r } \enspace,
\end{displaymath}
where $\lceil z \rceil$ denotes the smallest integer $N$ such that $z \leq N$.
Then we have $1 \leq p < \infty$ and $0 \leq \alpha < 1$ since $2/5 < r < 1$.
Let $\alpha = (0.\alpha_1 \alpha_2 \cdots)_2$ be the unique dyadic expansion of $\alpha$ with infinitely many zeroes.
By using these notations, we define finite binary words $w^{\langle 0 \rangle},w^{\langle 1 \rangle},\dots$ such that $w^{\langle i \rangle}$ is a proper initial subword of $w^{\langle i+1 \rangle}$ for each $i \geq 0$ by
\begin{displaymath}
w^{\langle 0 \rangle} = \emptyset \mbox{ and }
w^{\langle i \rangle} = w^{\langle i-1 \rangle}w^{\langle i-1 \rangle}010110^{p - \alpha_i} \mbox{ for } i \geq 1 \enspace.
\end{displaymath}
By the construction, the sequence $w^{\langle 0 \rangle},w^{\langle 1 \rangle},\dots$ converges to an infinite word.
We define $w$ to be the limit of the sequence.
We show that this $w$ satisfies the above condition.
For simplicity, put
\begin{displaymath}
\ell_i = \ell(w^{\langle i \rangle}) \enspace,\enspace \zeta_i = |Z(w^{\langle i \rangle})| \enspace,\enspace \pi_i = |P(w^{\langle i \rangle})| \mbox{ for each } i \geq 1 \enspace.
\end{displaymath}
These quantities are calculated as follows:
\begin{lemma}
\label{lem:tight_each_part}
For any $i \geq 1$, we have
\begin{displaymath}
\begin{split}
&\ell_i = (2^i - 1)(p + 5) - \sum_{j = 1}^{i} 2^{i-j} \alpha_j \enspace,\enspace
\zeta_i = (2^i - 1)(p + 2) - \sum_{j = 1}^{i} 2^{i-j} \alpha_j \enspace, \\
&\pi_i = (2^i - 1)p - 1 - \sum_{j = 1}^{i} 2^{i-j} \alpha_j + \delta_{p,1} \alpha_i \enspace,
\end{split}
\end{displaymath}
where $\delta_{a,b}$ denotes the Kronecker delta.
\end{lemma}
\begin{proof}
We prove the lemma by induction on $i$.
Since $w^{\langle 1 \rangle} = 010110^{p - \alpha_1}$, the case $i = 1$ is derived by a direct calculation.
From now, we consider the case $i \geq 2$.

For $\ell_i$ and $\zeta_i$, the construction of $w^{\langle i \rangle}$ and the induction hypothesis imply that
\begin{displaymath}
\begin{split}
\ell_i
&= 2 \ell_{i-1} + p - \alpha_i + 5 \\
&= (2(2^{i-1} - 1) + 1)(p + 5) - 2\sum_{j = 1}^{i-1} 2^{i-1-j} \alpha_j - \alpha_i
= (2^i - 1)(p + 5) - \sum_{j = 1}^{i} 2^{i-j} \alpha_j \enspace, \\
\zeta_i
&= 2 \zeta_{i-1} + p - \alpha_i + 2 \\
&= (2(2^{i-1} - 1) + 1)(p + 2) - 2\sum_{j = 1}^{i-1} 2^{i-1-j} \alpha_j - \alpha_i
= (2^i - 1)(p + 2) - \sum_{j = 1}^{i} 2^{i-j} \alpha_j \enspace,
\end{split}
\end{displaymath}
therefore the claim holds.
We deal with the value $\pi_i$ in the rest of this proof.
First, note that by the construction of $w^{\langle i \rangle}$, each one in $w^{\langle i \rangle}$ has one of the three properties; it is followed by a one, it is followed by $011$, and it is preceded by a one.
This implies that an index $j$ with $(w^{\langle i \rangle})_j = 1$ cannot be a member of the set $P(w^{\langle i \rangle})$.
Thus by the definition of $P(w^{\langle i \rangle})$, we have $j \in P(w^{\langle i \rangle})$ if and only if $j \geq 2$ and $(w^{\langle i \rangle})_{j-1} = (w^{\langle i \rangle})_j = 0$.

First we consider the case $p \geq 2$.
In this case, $w^{\langle i-1 \rangle}$ ends with a zero regardless of the value $\alpha_{i-1} \in \{0,1\}$, while $w^{\langle i-1 \rangle}$ begins with a zero.
Thus the members of $P(w^{\langle i \rangle})$ are classified into the four types; those induced by $P(w^{\langle i-1 \rangle})$ for the first two blocks $w^{\langle i-1 \rangle}$ in $w^{\langle i \rangle} = w^{\langle i-1 \rangle} w^{\langle i-1 \rangle} 010110^{p-\alpha_i}$, the first position in the second $w^{\langle i-1 \rangle}$, the first position in the last block $010110^{p - \alpha_i}$, and the last $p - \alpha_i - 1$ positions in the last block $010110^{p - \alpha_i}$.
Summarizing, we have
\begin{displaymath}
\begin{split}
\pi_i
&= 2 \pi_{i-1} + p - \alpha_i + 1 \\
&= (2(2^{i-1} - 1) + 1)p - 2 + 1 - 2\sum_{j = 1}^{i-1} 2^{i-1-j} \alpha_j - \alpha_i
= (2^i - 1)p - 1 - \sum_{j = 1}^{i} 2^{i-j} \alpha_j
\end{split}
\end{displaymath}
where we used the induction hypothesis to derive the second equality (recall that now $p \geq 2$), therefore the claim holds in this case.

Secondly, we consider the case $p = 1$.
In this case, $w^{\langle i-1 \rangle}$ ends with a zero if and only if $\alpha_{i-1} = 0$, while $w^{\langle i-1 \rangle}$ always begins with a zero.
Thus the members of $P(w^{\langle i \rangle})$ are classified into the three types; those induced by $P(w^{\langle i-1 \rangle})$ for the first two blocks $w^{\langle i-1 \rangle}$ in $w^{\langle i \rangle}$, the first position in the second $w^{\langle i-1 \rangle}$ (when $\alpha_{i-1} = 0$), and the first position in the last block $010110^{p - \alpha_i}$ (when $\alpha_{i-1} = 0$).
Summarizing, we have
\begin{displaymath}
\begin{split}
\pi_i = 2 \pi_{i-1} + 2 (1 - \alpha_{i-1}) 
&= 2 \left( 2^{i-1} - 2 - \sum_{j = 1}^{i-1} 2^{i-1-j} \alpha_j + \alpha_{i-1} \right) + 2 - 2 \alpha_{i-1} \\
&= 2^i - 2 - \sum_{j = 1}^{i-1} 2^{i-j} \alpha_j
= 2^i - 2 - \sum_{j = 1}^{i} 2^{i-j} \alpha_j + \alpha_i \enspace,
\end{split}
\end{displaymath}
where we used the induction hypothesis to derive the second equality (recall that now $p = 1$), therefore the claim holds in this case.
Hence Lemma \ref{lem:tight_each_part} holds.
\end{proof}
Lemma \ref{lem:tight_each_part} implies the following relation:
\begin{equation}
\label{eq:tight_relation_each_part}
5 \zeta_i - 2 \ell_i = 3 \pi_i + 3 - 3 \delta_{p,1} \alpha_i \mbox{ for any } i \geq 1 \enspace.
\end{equation}
Lemma \ref{lem:tight_each_part} gives the values of $|Z(w^{(n)})|$ and $|P(w^{(n)})|$ for special initial subwords $w^{(n)} = w^{\langle i \rangle}$ of $w$.
To consider a general case, we present the following decomposition of any (proper) initial subword of $w$:
\begin{lemma}
\label{lem:tight_decomposition}
Each finite initial subword $w^{(n)}$ of $w$ is decomposed as
\begin{equation}
\label{eq:tight_decomposition}
w^{(n)} = w^{\langle i_1 \rangle} w^{\langle i_2 \rangle} \cdots w^{\langle i_{k-1} \rangle} (w^{\langle i_k \rangle})^{\lambda+1} v \enspace,
\end{equation}
where $k \geq 1$, $i_1 > i_2 > \cdots > i_k \geq 0$, $\lambda \in \{0,1\}$, $v$ is a (possibly empty) initial subword of $010110^{p - \alpha_{i_k+1}}$, and $i_k \geq 1$ if $k \geq 2$.
\end{lemma}
\begin{proof}
Since $w$ is the limit of words $w^{\langle i \rangle}$, it suffices to show that every initial subword $u$ of each $w^{\langle m \rangle}$, $m \geq 0$, is expressed as in (\ref{eq:tight_decomposition}).
We proceed the proof by induction on $m$.
First, in the case $m \leq 1$, such an expression of $u$ is obtained by putting $k = 1$, $i_1 = 0$, $\lambda = 0$, and $v = u$.
From now, we focus on the case $m \geq 2$.

Recall that $w^{\langle m \rangle} = w^{\langle m-1 \rangle} w^{\langle m-1 \rangle} 010110^{p - \alpha_m}$ by the construction.
First, if $u$ is contained in the first $w^{\langle m-1 \rangle}$, then the claim follows from the induction hypothesis.
Secondly, if $u$ is not contained in the first two blocks $w^{\langle m-1 \rangle} w^{\langle m-1 \rangle}$, then we have $u = w^{\langle m-1 \rangle} w^{\langle m-1 \rangle} x$ for an initial subword $x$ of $010110^{p - \alpha_m}$, therefore a desired expression is obtained by putting $k = 1$, $i_k = m-1$, $\lambda = 1$ and $v = x$.
Finally, in the remaining case, we have $u = w^{\langle m-1 \rangle}x$ for a nonempty proper initial subword $x$ of $w^{\langle m-1 \rangle}$.
By the induction hypothesis, $x$ can be expressed in the following form:
\begin{displaymath}
x = w^{\langle i'_1 \rangle} w^{\langle i'_2 \rangle} \cdots w^{\langle i'_{k'-1} \rangle} (w^{\langle i'_{k'} \rangle})^{\lambda'+1} v' \enspace,
\end{displaymath}
where the parameters satisfy the same condition as in the statement of the lemma.
Note that $i'_1 \leq m-2$ by the choice of $x$.
Now a desired expression of $u$ is obtained by putting $k = k'+1$, $i_1 = m-1$, $i_j = i'_{j-1}$ for $2 \leq j \leq k$, $\lambda = \lambda'$, and $v = v'$.
Hence Lemma \ref{lem:tight_decomposition} holds.
\end{proof}
We use the decomposition (\ref{eq:tight_decomposition}) of $w^{(n)}$ in the following argument.
Moreover, we assume that $n \geq \ell_1$ since this is sufficient for our purpose.
Now we have $i_1 \geq 1$, therefore $i_k \geq 1$ (recall the condition that $i_k \geq 1$ if $k \geq 2$).
Thus we have
\begin{equation}
\label{eq:property_initial_subword}
n = \ell(w^{(n)}) = \sum_{j = 1}^{k} \ell_{i_j} + \lambda \ell_{i_k} + \ell(v) \enspace,\enspace
|Z(w^{(n)})| = \sum_{j = 1}^{k} \zeta_{i_j} + \lambda \zeta_{i_k} + |Z(v)| \enspace.
\end{equation}

For the value $|P(w^{(n)})|$, we consider the following two cases.
First suppose $p \geq 2$.
Then each $w^{\langle i_j \rangle}$ ends with a zero regardless of the value $\alpha_{i_j} \in \{0,1\}$, therefore the members of $P(w^{(n)})$ are classified into the six types; those induced by the sets $P(w^{\langle i_j \rangle})$ for $k + \lambda$ blocks $w^{\langle i_j \rangle}$, those induced by the set $P(v)$, the first positions in the $k - 2$ blocks $w^{\langle i_j \rangle}$ ($2 \leq j \leq k-1$), the first position in the first $w^{\langle i_k \rangle}$ (when $k \geq 2$), the first position in the second $w^{\langle i_k \rangle}$ (when $\lambda = 1$), and the first position in the last block $v$ (when $v \neq \emptyset$).
Summarizing, we have
\begin{displaymath}
|P(w^{(n)})| = \sum_{j = 1}^{k} \pi_{i_j} + \lambda \pi_{i_k} + |P(v)| + k + \lambda - \delta_{v,\emptyset} \mbox{ if } p \geq 2 \enspace.
\end{displaymath}

Secondly, suppose $p = 1$.
Then $w^{\langle i_j \rangle}$ ends with a zero if and only if $\alpha_{i_j} = 0$.
Now the members of $P(w^{(n)})$ are classified into the six types; those induced by the sets $P(w^{\langle i_j \rangle})$ for $k + \lambda$ blocks $w^{\langle i_j \rangle}$, those induced by the set $P(v)$, the first positions in the blocks $w^{\langle i_j \rangle}$ with $2 \leq j \leq k-1$ (when $\alpha_{i_{j-1}} = 0$), the first position in the first $w^{\langle i_k \rangle}$ (when $k \geq 2$ and $\alpha_{i_{k-1}} = 0$), the first position in the second $w^{\langle i_k \rangle}$ (when $\lambda = 1$ and $\alpha_{i_k} = 0$), and the first position in the last block $v$ (when $v \neq \emptyset$ and $\alpha_{i_k} = 0$).
This implies that
\begin{displaymath}
\begin{split}
|P(w^{(n)})|
& = \sum_{j = 1}^{k} \pi_{i_j} + \lambda \pi_{i_k} + |P(v)| + \sum_{j = 2}^{k} (1 - \alpha_{i_{j-1}}) + \lambda (1 - \alpha_{i_k}) + (1 - \delta_{v,\emptyset})(1 - \alpha_{i_k}) \\
& = \sum_{j = 1}^{k} \pi_{i_j} + \lambda \pi_{i_k} + |P(v)| + k - \sum_{j = 1}^{k} \alpha_{i_j} + (\lambda - \delta_{v,\emptyset})(1 - \alpha_{i_k}) \mbox{ if } p = 1 \enspace.
\end{split}
\end{displaymath}

Summarizing, regardless of the value $p \geq 1$, we have
\begin{equation}
\label{eq:property_initial_subword_2}
|P(w^{(n)})| = \sum_{j = 1}^{k} \pi_{i_j} + \lambda \pi_{i_k} + |P(v)| + k + (\lambda - \delta_{v,\emptyset})(1 - \delta_{p,1}\alpha_{i_k}) - \delta_{p,1} \sum_{j = 1}^{k} \alpha_{i_j} \enspace.
\end{equation}
By using the above results, we derive the following two properties:
\begin{lemma}
\label{lem:tight_ratio_0}
We have $\lim_{n \to \infty} |Z(w^{(n)})| / n = r$.
\end{lemma}
\begin{proof}
Let $n \geq \ell_1$.
First, since $\alpha - 2^{-i} < \sum_{j=1}^{i} 2^{-j} \alpha_j \leq \alpha$, Lemma \ref{lem:tight_each_part} implies that
\begin{equation}
\label{eq:bound_ell_zeta}
\begin{split}
&2^i (p + 5 - \alpha) - (p + 5) \leq \ell_i < 2^i (p + 5 - \alpha) - (p + 4) \enspace, \\
&2^i (p + 2 - \alpha) - (p + 2) \leq \zeta_i < 2^i (p + 2 - \alpha) - (p + 1) \enspace.
\end{split}
\end{equation}
Thus, by putting $A = \sum_{j = 1}^{k} 2^{i_j} + \lambda 2^{i_k}$, it follows from (\ref{eq:property_initial_subword}) that
\begin{displaymath}
\begin{split}
n &< \sum_{j = 1}^{k} \bigl( 2^{i_j} (p + 5 - \alpha) - (p + 4) \bigr) + \lambda \bigl( 2^{i_k} (p + 5 - \alpha) - (p + 4) \bigr) + \ell(v) \\
&= (p + 5 - \alpha) A - (p + 4)(k + \lambda) + \ell(v) \enspace,
\end{split}
\end{displaymath}
\begin{displaymath}
\begin{split}
|Z(w^{(n)})|
&\geq \sum_{j = 1}^{k} \bigl( 2^{i_j} (p + 2 - \alpha) - (p + 2) \bigr) + \lambda \bigl( 2^{i_k} (p + 2 - \alpha) - (p + 2) \bigr) + |Z(v)| \\
&= (p + 2 - \alpha) A - (p + 2)(k + \lambda) + |Z(v)| \enspace,
\end{split}
\end{displaymath}
therefore
\begin{equation}
\label{eq:tight_ratio_0_lbound}
\frac{ |Z(w^{(n)})| }{ n } > \frac{ (p + 2 - \alpha) A - (p+2)(k + \lambda) + |Z(v)| }{ (p + 5 - \alpha) A - (p+4)(k + \lambda) + \ell(v) } \enspace.
\end{equation}
Now note that $0 \leq |Z(v)| \leq \ell(v) \leq p + 5$ by the choice of $v$, therefore the values $|Z(v)|$ and $\ell(v)$ are bounded.
On the other hand, we have $k \leq i_1 \leq \log_2 A$, while $i_1 \to \infty$ when $n \to \infty$, therefore we have $A \to \infty$ and $k/A \to 0$ when $n \to \infty$.
Thus the right-hand side of (\ref{eq:tight_ratio_0_lbound}) converges to $(p + 2 - \alpha)/(p + 5 - \alpha) = r$ when $n \to \infty$.

A similar argument also implies, by using (\ref{eq:property_initial_subword}) and (\ref{eq:bound_ell_zeta}), that
\begin{displaymath}
\frac{ |Z(w^{(n)})| }{ n } < \frac{ (p + 2 - \alpha) A - (p+1)(k + \lambda) + |Z(v)| }{ (p + 5 - \alpha) A - (p+5)(k + \lambda) + \ell(v) } \enspace,
\end{displaymath}
and the right-hand side converges to $(p + 2 - \alpha)/(p + 5 - \alpha) = r$ when $n \to \infty$.
Hence Lemma \ref{lem:tight_ratio_0} holds.
\end{proof}
\begin{lemma}
\label{lem:tight_ratio_E}
We have $\lim_{n \to \infty} |P(w^{(n)})| / n = (5r-2)/3$.
\end{lemma}
\begin{proof}
Let $n \geq \ell_1$.
Then by (\ref{eq:tight_relation_each_part}), (\ref{eq:property_initial_subword}), and (\ref{eq:property_initial_subword_2}), we have
\begin{displaymath}
\begin{split}
5 |Z(w^{(n)})| - 2 n
&= \sum_{j = 1}^{k} (5 \zeta_{i_j} - 2 \ell_{i_j}) + \lambda (5 \zeta_{i_k} - 2 \ell_{i_k}) + 5 |Z(v)| - 2 \ell(v) \\
&= \sum_{j = 1}^{k} (3 \pi_{i_j} + 3 - 3 \delta_{p,1} \alpha_{i_j}) + \lambda (3 \pi_{i_k} + 3 - 3 \delta_{p,1} \alpha_{i_k}) + 5 |Z(v)| - 2 \ell(v) \\
&= 3 \left( \sum_{j = 1}^{k} \pi_{i_j} + \lambda \pi_{i_k} + k + \lambda (1 - \delta_{p,1} \alpha_{i_k}) - \delta_{p,1} \sum_{j = 1}^{k} \alpha_{i_j} \right) + 5 |Z(v)| - 2 \ell(v) \\
&= 3 \bigl( |P(w^{(n)})| - |P(v)| + \delta_{v,\emptyset}(1 - \delta_{p,1}\alpha_{i_k}) \bigr) + 5|Z(v)| - 2\ell(v) \enspace.
\end{split}
\end{displaymath}
Thus, by putting $B = 5 |Z(v)| - 2 \ell(v) - 3 |P(v)| + 3 \delta_{v,\emptyset} (1 - \delta_{p,1}\alpha_{i_k})$, we have
\begin{displaymath}
\frac{ |P(w^{(n)})| }{ n } = \frac{ 5 |Z(w^{(n)})| / n - 2 }{ 3 } - \frac{ B }{ 3n } \enspace.
\end{displaymath}
Note that $B$ is bounded since $\ell(v) \leq p + 5$.
Thus Lemma \ref{lem:tight_ratio_0} implies that the right-hand side converges to $(5r-2)/3$ when $n \to \infty$.
Hence Lemma \ref{lem:tight_ratio_E} holds.
\end{proof}
Hence, by Lemmas \ref{lem:tight_ratio_0} and \ref{lem:tight_ratio_E}, the infinite word $w$ satisfies the desired condition.

\subsection{Construction of Words Achieving Upper Bounds}
\label{subsec:proof_upper_bound_best}

We show that for any $0 \leq r \leq 1$, there exists an infinite binary word $w$ such that $r_{\inf}(w) = r_{\sup}(w) = r$ and $R_{\inf}(w) = R_{\sup}(w) = r$.
Note that $w = 0000 \cdots$ (infinite repetition of $0$) satisfies the condition for the case $r = 1$.
From now, we focus on the remaining case $0 \leq r < 1$.

To construct such a word $w$, first define $\delta_k \in \{0,1\}$ for $k \geq 1$ inductively by
\begin{displaymath}
\delta_k = 1\ \mbox{ if } \frac{ \sum_{i = 1}^{k-1} 2i \delta_i + 2k }{ k(k+1) } \leq r \enspace,\enspace \delta_k = 0\ \mbox{ otherwise.}
\end{displaymath}
Then it follows by the induction on $k$ that
\begin{equation}
\label{eq:upper_bestpossible_1}
\frac{ \sum_{i = 1}^{k} 2i\delta_i }{ k(k+1) } \leq r \mbox{ for any } k \geq 1 \enspace.
\end{equation}
Indeed, if $\delta_k = 1$, then the inequality holds by the definition of $\delta_k$; while if $\delta_k = 0$, then the induction hypothesis implies that
\begin{displaymath}
\frac{ \sum_{i=1}^{k} 2i\delta_i }{ k(k+1) } \leq \frac{ \sum_{i=1}^{k-1} 2i\delta_i }{ (k-1)k } \leq r
\end{displaymath}
(note that $\delta_1 = 0$ and the claim indeed holds when $k = 1$).
Now let $w^{\langle k \rangle}$ denote the repetition of $1 - \delta_k \in \{0,1\}$ of length $2k$ for every $k \geq 1$, and define $w = w^{\langle 1 \rangle} w^{\langle 2 \rangle} w^{\langle 3 \rangle} \cdots$.
We prove that this infinite word $w$ satisfies the desired condition.
Note that
\begin{equation}
\label{eq:upper_bound_word_property}
\ell(w^{\langle 1 \rangle} \cdots w^{\langle k \rangle}) = \sum_{i = 1}^{k} 2i = k(k+1)\ \mbox{ and }
|Z(w^{\langle 1 \rangle} \cdots w^{\langle k \rangle})| = \sum_{i = 1}^{k} 2i\delta_i \enspace.
\end{equation}
Now we have the following two properties:
\begin{lemma}
\label{lem:upper_bestpossible_1}
We have $\lim_{n \to \infty} |Z(w^{(n)})| / n = r$.
\end{lemma}
\begin{proof}
For an arbitrary $\varepsilon > 0$, take an integer $M > 0$ such that $1/M \leq \varepsilon / 2$.
Then by (\ref{eq:upper_bound_word_property}), we have
\begin{displaymath}
\lim_{n \to \infty} \frac{ \sum_{i = 1}^{M} 2i\delta_i + \sum_{i = M+1}^{n} 2i }{ n(n+1) } = \lim_{n \to \infty} \left( 1 - \frac{ \sum_{i = 1}^{M} 2i(1-\delta_i) }{ n(n+1) } \right) = 1 \enspace,
\end{displaymath}
while $r < 1$.
Thus by the definition of $\delta_k$, we have $\delta_K = 0$ for some $K > M$.
We show that $r - \varepsilon \leq |Z(w^{(n)})| / n \leq r$ for any $n \geq (K-1)K = \ell(w^{\langle 1 \rangle} \cdots w^{\langle K-1 \rangle})$, from which the claim of the lemma follows.

For the purpose, we show, by induction on $k \geq K$, that $r - \varepsilon \leq |Z(w^{(n)})| / n \leq r$ for any $(k-1)k \leq n \leq k(k+1)$.
Take such an $n$.
Then we have $w^{(n)} = w^{\langle 1 \rangle} \cdots w^{\langle k-1 \rangle} v$, where $v$ is a subword of $w^{\langle k \rangle}$ of length $n - (k-1)k$.
Now if $\delta_k = 0$ (that is satisfied when $k = K$ by the choice of $K$), then (\ref{eq:upper_bestpossible_1}) and (\ref{eq:upper_bound_word_property}) imply that
\begin{displaymath}
\frac{ |Z(w^{(n)})| }{ n }
= \frac{ |Z(w^{\langle 1 \rangle} \cdots w^{\langle k-1 \rangle})| }{ n }
= \frac{ \sum_{i = 1}^{k-1} 2i\delta_i }{ n }
\leq \frac{ \sum_{i = 1}^{k-1} 2i\delta_i }{ (k-1)k } \leq r \enspace,
\end{displaymath}
while the definition of $\delta_k$ and the fact $k \geq K > M$ imply that
\begin{displaymath}
\frac{ |Z(w^{(n)})| }{ n }
= \frac{ \sum_{i = 1}^{k-1} 2i\delta_i }{ n }
\geq \frac{ \sum_{i = 1}^{k-1} 2i\delta_i }{ k(k+1) }
> r - \frac{ 2 }{ k+1 } > r - \frac{ 2 }{ M } \geq r - \varepsilon
\end{displaymath}
where the last inequality follows from the choice of $M$.
Thus $r - \varepsilon \leq |Z(w^{(n)})| / n \leq r$ if $\delta_k = 0$.
On the other hand, if $\delta_k = 1$ (hence $k \geq K + 1$), then we have
\begin{displaymath}
\frac{ |Z(w^{(n)})| }{ n }
= \frac{ |Z(w^{\langle 1 \rangle} \cdots w^{\langle k-1 \rangle})| + \ell(v) }{ \ell(w^{\langle 1 \rangle} \cdots w^{\langle k-1 \rangle}) + \ell(v)}
\leq \frac{ |Z(w^{\langle 1 \rangle} \cdots w^{\langle k-1 \rangle})| + 2k }{ \ell(w^{\langle 1 \rangle} \cdots w^{\langle k-1 \rangle}) + 2k}
\end{displaymath}
where the last inequality holds since $|Z(w^{\langle 1 \rangle} \cdots w^{\langle k-1 \rangle})| \leq \ell(w^{\langle 1 \rangle} \cdots w^{\langle k-1 \rangle})$ and $0 \leq \ell(v) \leq 2k$.
Now it follows from (\ref{eq:upper_bestpossible_1}) and (\ref{eq:upper_bound_word_property}) that
\begin{displaymath}
\frac{ |Z(w^{(n)})| }{ n }
\leq \frac{ \sum_{i = 1}^{k-1} 2i\delta_i + 2k }{ (k-1)k + 2k }
= \frac{ \sum_{i = 1}^{k} 2i\delta_i }{ k(k+1) }
\leq r \enspace.
\end{displaymath}
Similarly, we have
\begin{displaymath}
\frac{ |Z(w^{(n)})| }{ n }
= \frac{ |Z(w^{\langle 1 \rangle} \cdots w^{\langle k-1 \rangle})| + \ell(v) }{ \ell(w^{\langle 1 \rangle} \cdots w^{\langle k-1 \rangle}) + \ell(v)}
\geq \frac{ |Z(w^{\langle 1 \rangle} \cdots w^{\langle k-1 \rangle})| }{ \ell(w^{\langle 1 \rangle} \cdots w^{\langle k-1 \rangle}) }
\end{displaymath}
by the properties $\ell(v) \geq 0$ and $|Z(w^{\langle 1 \rangle} \cdots w^{\langle k-1 \rangle})| \leq \ell(w^{\langle 1 \rangle} \cdots w^{\langle k-1 \rangle})$.
Thus, by putting $k' = k-1$, it follows from (\ref{eq:upper_bestpossible_1}) and (\ref{eq:upper_bound_word_property}) that
\begin{displaymath}
\frac{ |Z(w^{(n)})| }{ n }
\geq \frac{ \sum_{i = 1}^{k-1} 2i\delta_i }{ (k-1)k }
= \frac{ |Z(w^{(k'(k'+1))})| }{ k'(k'+1) }
\geq r - \varepsilon
\end{displaymath}
where the last inequality follows from the induction hypothesis.
Thus we have $r - \varepsilon \leq |Z(w^{(n)})| / n \leq r$ if $\delta_k = 1$.
Summarizing, we have $r - \varepsilon \leq |Z(w^{(n)})| / n \leq r$ regardless of the value $\delta_k$.
Hence Lemma \ref{lem:upper_bestpossible_1} holds.
\end{proof}
\begin{lemma}
\label{lem:upper_bestpossible_2}
We have $\lim_{n \to \infty} |P(w^{(n)})| / n = r$.
\end{lemma}
\begin{proof}
First, note that $w^{(n)}$ can be written as $w^{(n)} = w^{\langle 1 \rangle} \cdots w^{\langle k-1 \rangle} v$, where $v$ is a subword of $w^{\langle k \rangle}$ of length $n' = n - (k-1)k > 0$.
Now if $\delta_k = 0$, then we have $|Z(w^{(n)})| = \sum_{i = 1}^{k-1} 2i\delta_i$ by (\ref{eq:upper_bound_word_property}), while we have $|P(w^{(n)})| \geq \sum_{i = 1}^{k-1} (2i-1)\delta_i$ by counting the subwords $00$ in $w^{(n)}$.
On the other hand, if $\delta_k = 1$, then we have $|Z(w^{(n)})| = \sum_{i = 1}^{k-1} 2i\delta_i + n'$ by (\ref{eq:upper_bound_word_property}), while we have $|P(w^{(n)})| \geq \sum_{i = 1}^{k-1} (2i-1)\delta_i + n' - 1$ by counting the subwords $00$ in $w^{(n)}$.
Summarizing, regardless of the value $\delta_k$, we have
\begin{displaymath}
0 \leq \frac{ |Z(w^{(n)})| }{ n } - \frac{ |P(w^{(n)})| }{ n }
\leq \frac{ \sum_{i = 1}^{k-1} \delta_i + 1 }{ n } \leq \frac{ k }{ n }
\end{displaymath}
where the first inequality follows from the fact that there exists an injection $P(w^{(n)}) \to Z(w^{(n)})$ (see the first paragraph of Section \ref{sec:proof}).
Moreover, we have $k/n \to 0$ when $n \to \infty$, since $n > (k-1)k$ by the choice of the expression of $w^{(n)}$.
Thus we have
\begin{displaymath}
\lim_{n \to \infty} \frac{ |P(w^{(n)})| }{ n } = \lim_{n \to \infty} \frac{ |Z(w^{(n)})| }{ n } = r \enspace,
\end{displaymath}
where the second equality follows from Lemma \ref{lem:upper_bestpossible_1}.
Hence Lemma \ref{lem:upper_bestpossible_2} holds.
\end{proof}
By Lemmas \ref{lem:upper_bestpossible_1} and \ref{lem:upper_bestpossible_2}, the infinite word $w$ satisfies the desired condition.
Hence the proof of Theorem \ref{thm:main} is concluded.

\end{document}